\documentclass[11pt]{article}
\usepackage[utf8]{inputenc}
\usepackage{amsmath, amssymb, amsfonts, amsthm}
\usepackage[none]{hyphenat}
\usepackage{float}
\usepackage{mathtools}
\usepackage[colorlinks=true]{hyperref}
\hypersetup{
  colorlinks=true,
  linkcolor=blue,
  filecolor=magenta,
  urlcolor=cyan,
  citecolor=blue
}
\usepackage{cleveref}
\usepackage[top=25mm, bottom=25mm, left=25mm, right=25mm]{geometry}
\usepackage{xcolor}
\usepackage{comment}

\theoremstyle{plain}
\newtheorem{theorem}{Theorem}
\newtheorem{lemma}{Lemma}

\newtheorem{remark}{Remark}
\newtheorem{problem}{Problem}

\newcommand{\F}{\mathbb{F}}

\newcommand{\dist}{d_H}

\title{Triangle-free subsets of the $r$-distance graph of the hypercube}
\author{Padmini Mukkamala\thanks{BME Budapest Institute of Technology and Economics, Budapest, Hungary. Supported by staff mobility funded by ERASMUS+ EU fund}
\and
Ananthakrishnan Ravi\thanks{Delft University of Technology, Netherlands. Supported by an NWO open competition
grant (OCENW.M.22.090). \emph{E-mail}: \text{a.ravi@tudelft.nl}.}
}
\date{}

\begin{document}
\maketitle

\begin{abstract}
Given the $r$-distance graph on the hypercube $\F_2^n$, where two vertices are adjacent if their Hamming distance is exactly $r$, we study the maximum size $T(n,r)$ of a triangle-free set of vertices.
For even $r\le n/2$, we prove
\[
T(n,r)=O\!\left(\frac{r2^n}{n+1}\right).
\]
In particular, $T(n,r)=o(2^n)$ whenever $r=o(n)$.
For fixed $0<\alpha<2/3$, we also prove that if $r=\alpha n$, where $n$ ranges over integers such that $\alpha n$ is an even integer, then
\[
T(n,r)\le 2^{(1-\varepsilon_\alpha)n}
\]
for some $\varepsilon_\alpha>0$.
We also obtain lower bounds in various regimes of $r$ as a function of $n$.
\end{abstract}

\section{Introduction}\label{sec:intro}

We work in the setting of the hypercube $\F_2^n$.
For two vectors $x,y \in \F_2^n$, the Hamming distance between $x$ and $y$ is the number of coordinates at which they differ.
We denote this by $\dist(x,y)$. Fix an integer $r$ with $1\le r\le n$.
The \emph{$r$-distance graph} $H(n,r)$ has vertex set $\F_2^n$, and we join $x$ and $y$ by an edge when $\dist(x,y)=r$.
A \emph{triangle} in $H(n,r)$ is a triple $\{x,y,z\}$ with
\[
\dist(x,y)=\dist(x,z)=\dist(y,z)=r.
\]
Our goal is to understand how large a subset of $\F_2^n$ can be that does not contain triangles.
Let $T(n,r)$ be the maximum size of a triangle-free set of vertices in $H(n,r)$.

Two constraints appear immediately.
Triangles can occur only when $r$ is even, and only when $3r\le 2n$.
For any triple $x,y,z\in\F_2^n$, the sum $\dist(x,y)+\dist(x,z)+\dist(y,z)$ is always even, since each coordinate contributes either $0$ (all bits equal) or $2$ (exactly one bit differs).
The same coordinate-wise count also gives $\dist(x,y)+\dist(x,z)+\dist(y,z)\le 2n$, so for a triangle this sum equals $3r$ and hence $3r\le 2n$.
So we focus on even $r$ with $r\le 2n/3$.

Before we move to the general setting, we look at two small cases.
For $n=3$ and $r=2$, the set $\{{000},\ {001},\ {010},\ {011}\}$
is a maximum triangle-free set in $H(3,2)$. Indeed, the even-weight and odd-weight vertices each induce a copy
of $K_4$ in $H(3,2)$. A triangle-free set contains at most two
vertices from each parity class, and hence has at most four
vertices in total.
For $n=4$ and $r=2$, the set
\[
\{{0000},\ {0001},\ {0101},\ {0111},\ {1000},\ {1010},\ {1110},\ {1111}\}
\]
is a maximum triangle-free set in $H(4,2)$. To see maximality, observe that each parity class induces
$K_8$ with a perfect matching removed: the only nonadjacent
distinct pairs within a parity class are antipodal pairs at
distance $4$. Any five vertices in such a component contain a
triangle, so a triangle-free set contains at most four vertices
of each parity. It therefore has at most eight vertices in total.
These examples highlight that triangle-free sets in $H(n,r)$ need not be independent: they may contain many edges of length $r$ while still avoiding triangles.

The problem may be viewed as a higher-order exact-distance problem on the hypercube.
If we forbid all pairs at distance $r$, we obtain the independence number of the graph $H(n,r)$, which is the classical forbidden-distance problem; in the linear regime this is the setting of the Frankl--R\"odl graph \cite{frankl,kauer}.
In coding-theoretic language, this asks for the largest binary code in which the distance $r$ never occurs between two distinct codewords \cite{enomoto}.
Questions of this type are part of the broader study of distance constraints in the hypercube; see, for example, \cite{kleitman}.
Here we allow pairs at distance $r$, but forbid triples of vertices that are pairwise at distance $r$.
Equivalently, if $\mathcal H(n,r)$ denotes the $3$-uniform hypergraph on $\F_2^n$ whose hyperedges are the triangles of $H(n,r)$, then $T(n,r)$ is the size of the largest independent set in $\mathcal H(n,r)$.
In other words,
\[
T(n,r)=\alpha(\mathcal H(n,r)),
\]
where $\alpha(\mathcal H)$ denotes the maximum size of a vertex set containing no hyperedge of $\mathcal H$.

This hypergraph viewpoint is also useful for upper bounds.
Lov\'asz \cite{lovasz_theta} introduced the theta number as an efficiently computable upper bound for the independence number of a graph, and Gr\"otschel, Lov\'asz, and Schrijver \cite{lovasz} developed the related notion of the theta body.
More recently, Castro-Silva, de Oliveira Filho, Slot, and Vallentin \cite{valentin} introduced a recursive theta body for hypergraphs and applied it to this setting, obtaining an upper bound for $T(n,r)$.
Their results suggest that in the linear regime $r\sim n/c$, the behaviour should change around the threshold $r=n/2$: for $c>2$ the normalized upper bound is conjectured to decay exponentially fast, while at $c=2$ only linear decay is expected \cite[Conjecture~7.3]{valentin}.

The discussion above points toward upper-bound methods.
However, our first contribution is to give lower bounds showing how large triangle-free sets can be in different
regimes of $r$.
Since the triangle-free condition is weaker than independence, we can exploit constructions that allow many edges of length $r$ while still avoiding triangles.
This is reflected in our lower bounds, which come from three different mechanisms.
The first mechanism constructs independent sets by forbidding distance $r$, using minimum-distance codes.
The second keeps a large set and then deletes a small number of vertices to destroy every triangle, via the method of alteration.
The third restricts to a structured slice of the hypercube and enforces a local obstruction that prevents triangles from occurring.

We combine these three ideas to obtain lower bounds in several regimes of $r$.

\begin{theorem}\label{thm:lower-main}
Let $r$ be even and satisfy $r\le 2n/3$.
Then $T(n,r)$ satisfies the following lower bounds.

\begin{enumerate}
\item Fixed distance $r=2d$.
For every fixed integer $d\ge 1$,
\[
T(n,2d)=\Omega_d\!\left(\frac{2^n}{n^d}\right).
\]

\item If $r=o(n^{1/3})$ with $r\to\infty$, then
\[
T(n,r)\ \ge\ 2^{\,n-\frac r2\log_2 n-O(r)}.
\]

\item If $r=o(n)$ with $r\to\infty$, then
\[
T(n,r)\ \ge\ 2^{\,n-\frac{3}{4}r\log_2(n/r)-O(r)}.
\]

We further divide the case $r=\alpha n$ (where $n$ ranges over integers for which $r=\alpha n$ is an even integer) into two intervals, one where $\alpha\in(0,\alpha_0)$ and the other where $\alpha\in[\alpha_0,2/3]$.
Numerically, $\alpha_0\approx 0.16241686$.

\item If $r=\alpha n$ with fixed $0<\alpha<\alpha_0$, then
\[
T(n,\alpha n)\ \ge\ 2^{(c(\alpha)-o(1))n},
\]
where
\[
c(\alpha)=
1-\frac12\left(
H_2(\alpha)+\alpha+(1-\alpha)H_2\!\left(\frac{\alpha}{2(1-\alpha)}\right)
\right).
\]

\item If $r=\alpha n$ with fixed $\alpha_0\le \alpha\le 2/3$, then
\[
T(n,\alpha n)\ \ge\ 2^{(H_2(\alpha/2)-o(1))n},
\]
where
\[
H_2(t)=-t\log_2 t-(1-t)\log_2(1-t)
\]
is the binary entropy function.
\end{enumerate}

Parts (1)--(2) follow from \Cref{lem:bch}.
Parts (3)--(4) are asymptotic corollaries of the alteration bound proved in \Cref{lem:alteration}.
Part (5) follows from the two-bit construction in \Cref{lem:twobit}.
Here, $\alpha_0\in(0,2/3)$ is the threshold at which the exponents of \Cref{lem:alteration} and \Cref{lem:twobit} agree.
\end{theorem}

We also note that when $r$ is an even integer and $r=\frac{2n}{3}$, then the cube $\F_2^{n-1}$ does not contain any triangles, so there is a triangle-free subset of size $2^n/2$.
This construction extends near the endpoint $r=2n/3$.
More generally, if $r$ is an even integer and $r=\frac{2}{3}(n-c)$ for some integer $c\ge 0$, then any cube of codimension $c+1$ is triangle-free, since its dimension is $n-c-1$ and $
3r=2(n-c)>2(n-c-1)$.
Hence, $T(n,r)\ge 2^{n-c-1}$.

In particular, when $c=o(\log n)$ this gives
\[
T(n,r)\ge 2^{n-o(\log n)}=\frac{2^n}{n^{o(1)}},
\]
which is stronger than the lower bounds coming from our other constructions in this regime.

For upper bounds, we use two strategies.
The first uses the layered structure of the hypercube.
On each level, suitable shadows turn triangle-free property into a matching-type restriction, and a theorem of Frankl on families with no three pairwise disjoint sets bounds each level.

The second strategy is for the linear regime.
Fixing $x\in A$, the distance-$r$ neighbors of $x$ inside a triangle-free set $A$ are indexed by an $r$-uniform family avoiding intersection size $r/2$.
The Frankl--R\"odl forbidden-intersection theorem \cite{frankl} gives an exponentially small upper bound for each such local family.
Summing over $x\in A$ gives an upper bound on the number of distance-$r$ pairs in $A$, while a supersaturation theorem of Keevash and Long \cite{KeevashLong} gives a competing lower bound whenever $A$ is too large.
Comparing the two bounds gives the second upper bound.

\begin{theorem}\label{thm:upper-main}
Let $r$ be even.
Then $T(n,r)$ satisfies the following upper bounds.

\begin{enumerate}
\item If $r\le n/2$, then
\[
T(n,r)\ =\ O\!\left(\frac{r2^n}{n+1}\right).
\]
In particular, if $r=o(n)$, then
\[
T(n,r)=o(2^n).
\]

\item For every fixed $\alpha\in(0,2/3)$, there is a constant $\varepsilon_\alpha>0$ such that for $r=\alpha n$, 
\[
T(n,r)\ \le\ 2^{(1-\varepsilon_\alpha)n}
\]
for all sufficiently large $n$ such that $\alpha n$ is an even integer.
\end{enumerate}
\end{theorem}

\begin{remark}
After the first version of this paper \cite{mukk} appeared on arXiv,
Balogh, Chen, and Li \cite{balogh} proved that for every fixed even
$r=2d$,
\[
T(n,2d)=\Theta_d\!\left(\frac{2^n}{n^d}\right).
\]
Thus the fixed-distance lower bound in \Cref{thm:lower-main} has the
correct order of magnitude.  Their result also gives the sharp upper
bound when $r$ is fixed.  Therefore, \Cref{thm:upper-main}(1) is intended
to give a uniform bound for growing $r$, and in particular to show that
$T(n,r)=o(2^n)$ whenever $r=o(n)$.
\end{remark}

The rest of the paper is organized as follows.
Section~\ref{sec:lower} proves the lower bounds in \Cref{thm:lower-main}.
Section~\ref{sec:upper} proves \Cref{thm:upper-main}.
Section~\ref{sec:qary} discusses the $q$-ary setting and the larger-clique
version, and ends with open problems.

\section{Lower bounds}\label{sec:lower}

Throughout we fix even $r\le 2n/3$ and work in the graph $H(n,r)$.

\subsection{Minimum-distance codes}

An independent set in $H(n,r)$ is a set $C\subseteq \F_2^n$ such that
\[
\dist(x,y)\ne r
\qquad\text{for all distinct }x,y\in C.
\]
This set is automatically triangle-free.
Hence, any lower bound for the maximum size of a binary code avoiding the
single distance $r$ gives a lower bound for $T(n,r)$.

In this subsection we apply this idea to even distances $r=2d$.
We use the following stronger condition: instead of merely
forbidding the single distance $2d$, we use codes with minimum
distance at least $2d+1$. Such a code avoids every positive distance at
most $2d$, and therefore in particular avoids the distance $2d$.

Recall that a binary linear code with parameters $[n,k,\delta]$ is a
$k$-dimensional subspace of $\F_2^n$ whose minimum Hamming distance is at
least $\delta$.
For our purpose, it is enough to use the following standard lower bound
from classical BCH codes.

\begin{theorem}[{\cite[Ch.~9, \S1, Theorem~1, p.~258]{MS}}]\label{thm:bch}
For every integer $d\ge 1$, there is a binary BCH code of length
$N=2^m-1$ with dimension at least $N-dm$ and minimum distance at least
$2d+1$.
\end{theorem}

Since $m=\log_2(N+1)$, this gives a code of size at least
\[
2^{N-dm}
=
2^{N-d\log_2(N+1)}.
\]
We now transfer this bound to arbitrary lengths.

\begin{lemma}\label{lem:bch}
Let $d\ge 1$.
Then
\[
T(n,2d)\ \ge\ 2^{\,n-d\log_2(n+1)-d}.
\]
\end{lemma}

\begin{proof}
Choose $m=\lceil \log_2(n+1)\rceil$ and let $N=2^m-1$.
Then $N\ge n$ and $m\le \log_2(n+1)+1$.
By the BCH bound in \cref{thm:bch}, there is a binary linear code
$C\subseteq\F_2^N$ with
\[
\dim C\ge N-dm
\qquad\text{and}\qquad
d_{\min}(C)\ge 2d+1.
\]

Let $S\subseteq [N]$ have size $N-n$, and set
\[
C_0=\{c\in C:\ c_i=0 \text{ for every } i\in S\}.
\]
Then
\[
\dim C_0\ge \dim C-(N-n)\ge (N-dm)-(N-n)=n-dm.
\]
Let $C'$ be the code obtained from $C_0$ by deleting the coordinates in $S$. 
This does not change distances, since all these coordinates are zero in $C_0$. 
Thus, $C'\subseteq \F_2^n$ satisfies 
\[
|C'|\ge 2^{n-dm}
\qquad\text{and}\qquad
d_{\min}(C')\ge 2d+1.
\]

Hence $C'$ contains no pair of codewords at distance $2d$, and it is an independent set in $H(n,2d)$, therefore triangle-free.
Thus
\[
T(n,2d)\ge |C'|
\ge 2^{n-dm}
\ge 2^{n-d\log_2(n+1)-d}.
\]
\end{proof}

\subsection{A probabilistic bound}

We now prove the alteration bound.
We view triangles as hyperedges in a $3$-uniform hypergraph and remove one vertex from each surviving triangle.

\begin{lemma}\label{lem:alteration}
Let $r$ be even and satisfy $r\le 2n/3$.
Then $H(n,r)$ has a triangle-free vertex set of size at least
\[
\frac{2\sqrt{2}}{3}\cdot
\frac{2^n}{\sqrt{\binom{n}{r}\binom{r}{r/2}\binom{n-r}{r/2}}}.
\]
\end{lemma}

\begin{proof}
Let $\mathcal{H}(n,r)$ be the $3$-uniform hypergraph on the vertex set $\F_2^n$ whose hyperedges are the triangles of $H(n,r)$.
We choose each vertex independently with probability $p$.
Let $X$ be the number of chosen vertices and let $Y$ be the number of hyperedges of $\mathcal{H}(n,r)$ that are fully contained in the chosen set.

Removing one vertex from each such hyperedge produces a triangle-free set of size at least $X-Y$.
By linearity of expectation,
\[
\mathbb{E}[X-Y]\ =\ 2^n p - p^3\cdot |E(\mathcal{H}(n,r))|.
\]

It remains to count $|E(\mathcal{H}(n,r))|$.
Fix a vertex $u$.
There are $\binom{n}{r}$ vertices $v$ with $\dist(u,v)=r$.
Fix such a neighbor $v$.
A third vertex $w$ forms a triangle with $(u,v)$ exactly when $w$ differs from $u$ in $r/2$ of the $r$ coordinates where $u$ and $v$ differ, and also differs from $u$ in $r/2$ of the remaining $n-r$ coordinates.
So there are
\[
\binom{r}{r/2}\binom{n-r}{r/2}
\]
choices for $w$.
This counts each triangle $6$ times, so
\[
|E(\mathcal{H}(n,r))|
=\frac{1}{6}\,2^n\binom{n}{r}\binom{r}{r/2}\binom{n-r}{r/2}.
\]

Substituting this into $\mathbb{E}[X-Y]$ and optimizing over $p$ gives
\[
p=\sqrt{\frac{2}{\binom{n}{r}\binom{r}{r/2}\binom{n-r}{r/2}}}.
\]
With this choice,
\[
\mathbb{E}[X-Y]\ \ge\
\frac{2\sqrt{2}}{3}\cdot
\frac{2^n}{\sqrt{\binom{n}{r}\binom{r}{r/2}\binom{n-r}{r/2}}}.
\]
So there exists a choice of vertices with $X-Y$ at least this large, and the set obtained after deletion is triangle-free by construction.
\end{proof}

\subsection{A two-bit construction}

The next bound works best when $r$ is a constant fraction of $n$.
We place all vertices on the level $r/2$ and enforce a two-bit condition.

\begin{lemma}\label{lem:twobit}
Let $r$ be even.
Then $H(n,r)$ has a triangle-free vertex set of size at least
\[
\binom{n}{r/2}-\binom{n-2}{r/2}.
\]
\end{lemma}

\begin{proof}
Identify a vertex of $\F_2^n$ with the subset of $[n]$ given by its $1$-coordinates.
Work on the level $\binom{[n]}{r/2}$.
Let $\mathcal{F}$ be the family of all $(r/2)$-subsets that contain $1$ or contain $2$.
Equivalently, the first or the second bit is $1$.
Then
\[
|\mathcal{F}|=\binom{n}{r/2}-\binom{n-2}{r/2}.
\]

We claim that $\mathcal{F}$ is triangle-free in $H(n,r)$.
Suppose $A,B,C\in\mathcal{F}$ form a triangle.
Since all three sets have size $r/2$, the condition $\dist(A,B)=r$ is equivalent to $A\cap B=\emptyset$.
So a triangle in $H(n,r)$ inside this level would force $A,B,C$ to be pairwise disjoint.
That is impossible in $\mathcal{F}$, because among three sets each containing $1$ or $2$, two of them must share one of these two elements.
So $\mathcal{F}$ spans no triangle.
\end{proof}

\subsection{Proof of \Cref{thm:lower-main}}

\begin{proof}[Proof of \Cref{thm:lower-main}]
The five parts of \Cref{thm:lower-main} are proved by the three mechanisms established above: the minimum-distance code in \Cref{lem:bch}, the alteration bound in \Cref{lem:alteration}, and the two-bit construction in \Cref{lem:twobit}.
Hence $T(n,r)$ is at least the maximum of the relevant corresponding quantities.

\Cref{lem:bch} gives
\[
T(n,r)\ge 2^{\,n-\frac r2\lceil \log_2(n+1)\rceil}.
\]
For fixed $r=2d$, this implies
\[
T(n,2d)
\ge 2^{\,n-d\log_2(n+1)-d}
=
\Omega_d\!\left(\frac{2^n}{n^d}\right).
\]
This proves (1).

If $r=o(n^{1/3})$ and $r\to\infty$, then
\[
T(n,r)\ge 2^{\,n-\frac r2\lceil \log_2(n+1)\rceil}
       =2^{\,n-\frac r2\log_2 n-O(r)}.
\]
This proves (2).

We now analyze the alteration bound.
Using
\[
\binom{n}{r}\le \left(\frac{en}{r}\right)^r,
\qquad
\binom{r}{r/2}\le 2^r,
\qquad
\binom{n-r}{r/2}\le \left(\frac{2e(n-r)}{r}\right)^{r/2},
\]
\Cref{lem:alteration} gives, when $r=o(n)$ and $r\to\infty$,
\[
T(n,r)
\ge
2^{\,n-\frac{3}{4}r\log_2(n/r)-O(r)}.
\]
This proves (3).

When $r=\alpha n$ with fixed $0<\alpha<2/3$, we use
\[
\binom{n}{r}
\le
\frac{1}{\sqrt{2\pi n\alpha(1-\alpha)}}2^{nH_2(\alpha)}
\le
c\,2^{nH_2(\alpha)},
\]
where
\[
H_2(t)=-t\log_2 t-(1-t)\log_2(1-t)
\]
is the binary entropy function, and $c$ is a constant.
Then \Cref{lem:alteration} gives
\[
T(n,\alpha n)
\ge
2^{\left(
1-\frac12\left[
H_2(\alpha)+\alpha+(1-\alpha)H_2\!\left(\frac{\alpha}{2(1-\alpha)}\right)
\right]
-o(1)\right)n}.
\]
This proves (4).

Finally, \Cref{lem:twobit} gives
\[
T(n,r)
\ge
\binom{n}{r/2}\left(1-\frac{(n-r/2)^2}{n^2}\right).
\]
For $r=\alpha n$, we use
\[
\binom{n}{r/2}
\ge
\frac{2^{nH_2(\alpha/2)}}{n+1}.
\]
Therefore
\[
T(n,\alpha n)\ge 2^{(H_2(\alpha/2)-o(1))n}.
\]
This proves (5) whenever the two-bit exponent dominates.

To find the precise constant $\alpha_0$ which separates the two intervals where the above terms dominate, we define
\[
D(\alpha)
=
1-\frac12\left(
H_2(\alpha)+\alpha
+(1-\alpha)H_2\!\left(\frac{\alpha}{2(1-\alpha)}\right)
\right)
-H_2(\alpha/2).
\]
Putting $x=\alpha/2$, direct differentiation gives
\[
D'(\alpha)
=
\frac12\log_2\left(
\frac{x^5}{(1-3x)^3(1-x)^2}
\right).
\]
The function
\[
x\longmapsto
\frac{x^5}{(1-3x)^3(1-x)^2}
\]
is strictly increasing on $(0,1/3)$. Thus $D'$ changes sign at
most once. 
Moreover,
\[
\lim_{\alpha\to0^+}D(\alpha)=1>0,
\qquad
D(2/3)
=
1-\frac12\left(H_2(2/3)+\frac23\right)
-H_2(1/3)
<0.
\]
Since $D'$ changes sign at most once, from negative to positive,
the function $D$ first decreases and then increases. Its value at
the right endpoint is still negative, so it cannot cross zero on
the increasing part. It therefore has exactly one zero in
$(0,2/3)$.
We denote this zero by
$\alpha_0$. Numerically,
\[
\alpha_0\approx0.16241686.
\]
Consequently, the alteration exponent dominates for
$0<\alpha<\alpha_0$, whereas the two-bit exponent dominates for
$\alpha_0\le\alpha\le2/3$.

\end{proof}

\section{Upper bounds}\label{sec:upper}
We now prove the two upper bounds in \Cref{thm:upper-main}.
The first is the slice bound.
The second uses the Frankl--R\"odl theorem and the supersaturation theorem of Keevash and Long.

For our first proof, we introduce the following notions.

\medskip
\noindent\textbf{Levels, $m$-shadows, and $m$-covers.}
For $0\le k\le n$, let
\[
L_k:=\{x\in\F_2^n:\ |x|=k\}
\]
be the $k$th level of the hypercube, where $|x|$ denotes the Hamming weight of $x$ (the number of $1$'s in $x$).
Fix an integer $m\ge 1$.
Given $x\in L_k$, an \emph{$m$-shadow} of $x$ is obtained by changing exactly $m$ of its $1$'s to $0$'s.
Thus every $m$-shadow lies in $L_{k-m}$ and $x$ has exactly $\binom{k}{m}$ $m$-shadows.
Dually, an \emph{$m$-cover} of $x$ is obtained by changing exactly $m$ of its $0$'s to $1$'s.
Thus every $m$-cover lies in $L_{k+m}$ and $x$ has exactly $\binom{n-k}{m}$ $m$-covers.

We first prove the case $r=2$ as a warm-up, since it contains the shadow-counting idea.
\begin{theorem}\label{thm:upper-r2}
Any triangle-free subset of $H(n,2)$ has size at most $\frac{4\cdot 2^n}{n}$.
\end{theorem}

\begin{proof}
Let $A\subseteq \F_2^n$ be triangle-free in $H(n,2)$ and define $S_k:=|A\cap L_k|$.
We bound $S_k$ for each $k$ and then sum over $k$.

The level $L_0$ contributes one vertex, so we treat it separately.
Fix $k\ge 1$ and consider vertices in $A\cap L_k$.
Each vertex in $L_k$ has exactly $\binom{k}{1}=k$  $1$-shadows in $L_{k-1}$.
We claim that any vertex in $L_{k-1}$ can be a $1$-shadow of at most two vertices of $A\cap L_k$.
Indeed, if three distinct vertices in $L_k$ had the same $1$-shadow, then they pairwise differ only in
which one-bit was removed, so they have pairwise Hamming distance $2$ and form a triangle in $H(n,2)$.

Therefore, double-counting pairs $(x,y)$, where $x\in A\cap L_k$
and $y\in L_{k-1}$ is a $1$-shadow of $x$, gives
\[
kS_k \ \le\ 2\binom{n}{k-1},
\qquad\text{so}\qquad
S_k \ \le\ \frac{2}{k}\binom{n}{k-1}.
\]
Summing over $k$,
\[
|A|
\le 1+\sum_{k=1}^n \frac{2}{k}\binom{n}{k-1}
=1+\sum_{k=1}^n \frac{2}{n+1}\binom{n+1}{k}
\le 1+\frac{2}{n+1}\cdot 2^{n+1}
\le \frac{4\cdot 2^n}{n}.
\]
\end{proof}

In the above summation, we could have used $1$-shadows when $k\le n/2$ and then $1$-covers when $k>n/2$.
The analysis for covers is symmetric to the analysis for shadows, and gives the same final bound.
Although it is not necessary for $r=2$, we will prefer this viewpoint for the general upper bound.

Before we proceed to the bound for general $r$, we need a set-theoretic ingredient.
Identifying a vertex of the hypercube with its support, we may view vertices on a fixed level as $k$-subsets of $[n]$.
In the argument below, triangle-freeness will force certain such families to contain no three pairwise disjoint sets.
This places us in the setting of the Erd\H{o}s Matching problem \cite{erdos}, which asks for the largest $k$-uniform family with no $s$ pairwise disjoint sets.
We will use the following theorem due to Frankl \cite{frankl2}.

\begin{theorem}[{\cite[Theorem~10.3, p.~22]{frankl2}}]\label{thm:frankl-matching}
Let $s,t\ge2$ be integers. Suppose that $\mathcal{F}\subseteq \binom{X}{t}$, where $|X|\ge ts$, and $\mathcal{F}$ contains no $s$ pairwise disjoint sets.
Then
\[
|\mathcal{F}|\ \le\ (s-1)\binom{|X|-1}{t-1}.
\]
\end{theorem}

\begin{proof}[\textbf{Proof of \Cref{thm:upper-main}(1)}]

Let $A\subseteq \F_2^n$ be triangle-free in $H(n,r)$ and define
$S_k:=|A\cap L_k|$.
The case $r=2$ follows from \Cref{thm:upper-r2}, so we may assume that
$r\ge4$.
We bound $S_k$ for each $k$ and then sum over $k$.
By symmetry, it is enough to bound the levels $k\le n/2$ and
then double the resulting estimate.
For $k<r/2$, we use the trivial bound $S_k\le \binom nk$.

Now fix $k$ with $r/2\le k\le n/2$ and consider vertices in $A\cap L_k$.
For $x\in L_k$, consider its $\frac r2$-shadows in $L_{k-\frac r2}$.
Thus each $x\in L_k$ has exactly $\binom{k}{r/2}$ 
$\frac r2$-shadows.

Fix an $\frac r2$-shadow $y\in L_{k-\frac r2}$.
Any $x\in L_k$ having $\frac r2$-shadow $y$ is obtained by choosing an $\frac r2$-subset of the $0$-coordinates of $y$
and turning them into $1$'s.
Equivalently, after fixing $y$, the possible vertices $x$ correspond to $\frac r2$-subsets of a ground set of size
\[
n-\Bigl(k-\frac r2\Bigr).
\]
Among vertices of $A\cap L_k$ with the same $\frac r2$-shadow $y$, we cannot have three choices whose added $\frac r2$-subsets are pairwise disjoint, because those three vertices would then have pairwise Hamming distance $r$ and would form a triangle in $H(n,r)$.
Therefore, for each fixed $\frac r2$-shadow $y$, the family of added $\frac r2$-subsets contains no three pairwise disjoint sets.

We now apply \Cref{thm:frankl-matching} with $s=3$ and $t=\frac r2$.
This requires the ground set to have size at least $3\frac r2$, namely
\[
n-(k-\frac r2)\ \ge\ 3\frac r2.
\]
Under this hypothesis, \Cref{thm:frankl-matching} implies that for each fixed $\frac r2$-shadow $y$,
the number of vertices in $A\cap L_k$ having $\frac r2$-shadow $y$ is at most
\[
2\,\cdot\binom{n-(k-\frac r2)-1}{\frac r2-1}.
\]

Now double count pairs $(x,y)$ where $x\in A\cap L_k$ and $y$ is an $\frac r2$-shadow of $x$.
There are $\binom{k}{\frac r2}S_k$ such pairs.
On the other hand, there are $\binom{n}{k-\frac r2}$ possible $\frac r2$-shadows, and for each $\frac r2$-shadow at most
$2\cdot \binom{n-(k-\frac r2)-1}{\frac r2-1}$ choices of $x$.
Hence
\[
\binom{k}{\frac r2}S_k
\ \le\
2\binom{n}{k-\frac r2}\binom{n-(k-\frac r2)-1}{\frac r2-1}.
\]
Simplifying this inequality gives
\[
S_k\ \le\  \binom{n}{k}\cdot \frac{r}{n-(k-r/2)}.
\]

When $k\le n/2$ and $2r\le n$, we have
\[
n-(k-\frac r2)=n-k+\frac r2\ \ge\ \frac{n}{2}+\frac r2\ \ge\ 2\frac r2+\frac r2\ =\ 3\frac r2,
\]
so the hypothesis $n-(k-\frac r2)\ge 3\frac r2$ that is needed for \Cref{thm:frankl-matching} holds in this range.
When $k>n/2$, we repeat the same argument using $\frac r2$-covers instead of $\frac r2$-shadows, which gives the same
bound with $k$ replaced by $n-k$.
Therefore, when summing over all levels, we sum up to $\lfloor n/2 \rfloor$ and then double the result. When $n$ is even, this counts the middle level twice, which is harmless for an upper bound.

Using the trivial bound $S_k\le \binom{n}{k}$ for $k<\frac r2$, and the
shadow bound for $\frac r2\le k\le n/2$, we obtain
\begin{align*}
|A|
=\sum_{k=0}^n S_k
&\le
2\left(\sum_{k=0}^{\frac r2-1}\binom{n}{k}
+\sum_{k=\frac r2}^{\lfloor n/2 \rfloor}\binom{n}{k}\cdot\frac{r}{n-(k-r/2)}\right)\\
&\le
2\left(\sum_{k=0}^{\frac r2}\binom{n}{k}
+\sum_{k=\frac r2}^{\lfloor n/2 \rfloor }\binom{n}{k}\cdot\frac{r}{n-k+1}\right)\\
&\le
2\left(\sum_{k=0}^{\frac r2}\binom{n}{k}
+\frac{r}{n+1}\sum_{k=\frac r2}^{\lfloor n/2 \rfloor }\binom{n+1}{k}\right)=O\!\left(\frac{r2^n}{n+1}\right).
\end{align*}

We note that since $r\le \lfloor n/2 \rfloor$, we have $r/2\le \lfloor n/4 \rfloor$, and hence
\[
\sum_{k=0}^{r/2}\binom{n}{k}
\le
\sum_{k=0}^{\lfloor n/4\rfloor}\binom{n}{k}
\le
2^{H_2(1/4)n}
=
O\!\left(\frac{2^n}{n+1}\right).
\]
Therefore the contribution of the initial levels is absorbed
uniformly into the asserted bound.

\end{proof}

We now prove the upper bound in the linear regime.
We first state the two results that we need.
The first is the Frankl--R\"odl forbidden-intersection theorem, in the form stated as Theorem~1.1 in Keevash and Long \cite{KeevashLong}.
The second is the supersaturated version of the forbidden-distance theorem, stated as Theorem~1.10 in the same paper.

A family $\mathcal A$ of sets is called $l$-avoiding if
\[
|A\cap B|\ne l
\qquad\text{for all distinct }A,B\in\mathcal A.
\]

\begin{theorem}[Frankl--R\"odl, Theorem~1.1 of Keevash--Long]\label{thm:KL-FR}
Let $\alpha,\epsilon\in(0,1)$ with $\epsilon\le \alpha/2$.
Let $k=\lfloor \alpha n\rfloor$ and
\[
l\in [\max(0,2k-n)+\epsilon n,\ k-\epsilon n].
\]
Then any $l$-avoiding family $\mathcal A\subset \binom{[n]}{k}$ satisfies
\[
|\mathcal A|\le (1-\delta)^n\binom{n}{k}
\]
where $\delta=\delta(\alpha,\epsilon)>0$.
\end{theorem}

For the second result, let $[q]=\{1,\ldots,q\}$.
For $C\subseteq [q]^n$, let $d_H(x,y)$ denote the Hamming distance between $x$ and $y$.

\begin{theorem}[Keevash--Long, Theorem~1.10]\label{thm:KL-supersat}
Given $\tilde\epsilon,\eta\in(0,1/2)$ there is $\delta'>0$ such that the following holds.
Let $q\ge 2$ be an integer and let $C\subset [q]^n$ satisfy 
\[
|C|>q^{(1-\delta')n}.
\]
Suppose that $d\in\mathbb N$ satisfies
\[
\tilde\epsilon n<d<(1-\tilde\epsilon)n
\]
and that $d$ is even when $q=2$.
Then there are at least
\[
\binom{n}{d}(q-1)^d |C|q^{-\eta n}
\]
unordered pairs $x,y\in C$ with $d_H(x,y)=d$.
\end{theorem}

\begin{proof}[\textbf{Proof of \Cref{thm:upper-main}(2)}]
Fix $\alpha\in(0,2/3)$. We prove that there is
$\varepsilon_\alpha>0$ such that, for $r=\alpha n$, every triangle-free
set $A\subseteq\F_2^n$ satisfies
\[
|A|\le 2^{(1-\varepsilon_\alpha)n}
\]
for all sufficiently large $n$.

Let $A\subseteq\F_2^n$ be triangle-free in $H(n,r)$. For $x\in A$, define
\[
\mathcal F_x=\left\{S\in\binom{[n]}{r}:\ x+\mathbf 1_S\in A\right\}.
\]
Then $|\mathcal F_x|=|\{y\in A:\dist(x,y)=r\}|$. If
$S,T\in\mathcal F_x$, with corresponding vertices
$y=x+\mathbf 1_S$ and $z=x+\mathbf 1_T$, then
\[
\dist(y,z)=2r-2|S\cap T|.
\]
Thus $\dist(y,z)=r$ if and only if $|S\cap T|=r/2$. Since $A$ is
triangle-free, no two distinct sets in $\mathcal F_x$ have intersection
size $r/2$. Hence $\mathcal F_x$ is an $(r/2)$-avoiding family in
$\binom{[n]}{r}$.

Choose
\[
0<\epsilon<\min\left\{\frac{\alpha}{2},\ 1-\frac{3\alpha}{2}\right\},
\]
which is possible because $\alpha<2/3$. For all sufficiently large $n$,
we have
\[
\frac r2\in [\max(0,2r-n)+\epsilon n,\ r-\epsilon n].
\]
So \Cref{thm:KL-FR}, applied with $k=r$ and $l=r/2$, gives a constant
$\delta=\delta(\alpha,\epsilon)>0$ such that, for every $x\in A$,
\[
|\mathcal F_x|\le (1-\delta)^n\binom{n}{r}.
\]
Set
\[
\gamma=-\log_2(1-\delta)>0.
\]
Then
\[
|\mathcal F_x|
\le
2^{-\gamma n}\binom{n}{r}.
\]
Let
\[
P_r(A)=|\{(x,y)\in A^2:\dist(x,y)=r\}|.
\]
Summing over $x\in A$ gives
\[
P_r(A)=\sum_{x\in A}|\mathcal F_x|
\le |A|\,2^{-\gamma n}\binom{n}{r}.
\tag{1}
\]

We now use \Cref{thm:KL-supersat} for the lower bound. Choose
$0<\tilde\epsilon<\min\{\alpha,1-\alpha,1/2\}$ and $0<\eta<\min\{\gamma,1/2\}$.
Since $r=\alpha n$, we have
$\tilde\epsilon n<r<(1-\tilde\epsilon)n$ for all sufficiently large $n$.
Applying \Cref{thm:KL-supersat} with $q=2$, $C=A$, and $d=r$, there is a
constant $\delta'=\delta'(\tilde\epsilon,\eta)>0$ such that, if
$|A|>2^{(1-\delta')n}$, then
\[
P_r(A)\ge 2|A|\binom{n}{r}2^{-\eta n}.
\tag{2}
\]
Here the factor $2$ appears because $P_r(A)$ counts ordered pairs, whereas \Cref{thm:KL-supersat} counts unordered pairs.

Since $\eta<\gamma$, the lower bound in $(2)$ contradicts the upper bound
in $(1)$. Therefore $|A|>2^{(1-\delta')n}$ is impossible. Taking
$\varepsilon_\alpha=\delta'$ gives
\[
T(n,r)\le 2^{(1-\varepsilon_\alpha)n}.
\]
This proves \Cref{thm:upper-main}(2).
\end{proof}

\section{Discussion: the $q$-ary hypercube}\label{sec:qary}

Let $q\ge2$ be an integer, and write $[q]=\{0,1,\ldots,q-1\}$.
Let $H_q(n,r)$ be the graph on $[q]^n$ in which two vertices are
adjacent if their Hamming distance is exactly $r$. Let $T_q(n,r)$
be the maximum size of a subset of $[q]^n$ containing no triangle
in $H_q(n,r)$.

When we use linear codes over $\F_q$, we will additionally assume
that $q$ is a prime power and identify the alphabet with $\F_q$.

Triangles in $H_2(n,r)$ can occur only when $r$ is even and when
$3r\le 2n$.
Both restrictions are specific to the binary cube.
When $q\ge 3$, triangles exist in $H_q(n,r)$ for every $1\le r\le n$.

For $q=2$, the fixed even distance regime is now understood.
Balogh, Chen, and Li proved \cite{balogh} that for every fixed even
$r=2d$,
\[
T_2(n,2d)=\Theta\!\left(\frac{2^n}{n^d}\right).
\]
It is therefore natural to ask for the corresponding fixed-distance
behaviour in the $q$-ary cube.

We first record the case $r=1$.

\begin{lemma}
For every prime power $q\ge 3$,
\[
T_q(n,1)=2q^{n-1}.
\]
\end{lemma}

\begin{proof}
Let $u,v,w\in\F_q^n$.  
The three vertices form a triangle in
$H_q(n,1)$ if and only if there is a coordinate position $i\in[n]$ such
that
\[
u_j=v_j=w_j \qquad\text{for every } j\ne i,
\]
and $u_i,v_i,w_i$ are three distinct elements of $\F_q$.

Let $A\subseteq\F_q^n$ be triangle-free.  
Choose a coordinate position $i\in[n]$, and choose values $a_j\in\F_q$ for every $j\ne i$.
Consider the set of vertices
\[
\{x\in\F_q^n:\ x_j=a_j \text{ for every } j\ne i\}.
\]
The set $A$ contains at most two points from this set.  
Otherwise, $A$ would contain three distinct points from this set.
These three points agree in all coordinates except $i$.
Since they are distinct, their $i$th coordinates are also distinct.
Hence, they form a triangle in $H_q(n,1)$.

There are $q^{n-1}$ choices for the values $a_j$, $j\ne i$, so
\[
|A|\le 2q^{n-1}.
\]

For the lower bound, let
\[
A=\{x\in\F_q^n:\ x_1+\cdots+x_n\in\{0,1\}\}.
\]
Choose a coordinate position $i\in[n]$, and choose values $a_j\in\F_q$
for every $j\ne i$.  On the set of vertices with $x_j=a_j$ for every
$j\ne i$, the linear form $L(x)=x_1+\cdots+x_n$ has the form
\[
L(x)=x_i+c,
\]
where $c=\sum_{j\ne i}a_j$ is fixed.  As $x_i$ varies through $\F_q$,
the value $x_i+c$ also varies through all of $\F_q$ exactly once.
Hence exactly one point has $L(x)=0$, and exactly one point has
$L(x)=1$.

Thus $A$ contains exactly two points from every set obtained by fixing
all coordinates except the $i$th coordinate.  Since every triangle in
$H_q(n,1)$ must lie inside one such set, $A$ is triangle-free.
Finally, once $x_1,\ldots,x_{n-1}$ are fixed, there are exactly two
choices of $x_n$ for which $x_1+\cdots+x_n\in\{0,1\}$.

Hence
\[
|A|=2q^{n-1}.
\]
Therefore
\[
T_q(n,1)=2q^{n-1}.
\]
\end{proof}

Next we look at the $r=2$ case.

\begin{lemma}\label{lem:qary-r2}
For every fixed prime power $q\ge 3$,
\[
T_q(n,2)=\Theta_q\!\left(\frac{q^n}{n}\right).
\]
\end{lemma}
\begin{proof}
Let $A\subseteq\F_q^n$ be triangle-free in $H_q(n,2)$.
Fix $y\in\F_q^n$.  For each $i\in[n]$, let
\[
B_i(y)=\{x\in\F_q^n:\ x_j=y_j \text{ for every } j\ne i,\ x_i\ne y_i\}.
\]
The vertices at distance $1$ from $y$ are the disjoint union of the sets
$B_i(y)$, and each $B_i(y)$ has size $q-1$.

We claim that $A$ meets at most two of the sets $B_i(y)$.  Otherwise, if
$A$ contained points from three different sets $B_i(y),B_j(y),B_k(y)$,
then these three points would be pairwise at Hamming distance $2$, and so
would form a triangle in $H_q(n,2)$.  Hence, for every $y\in\F_q^n$,
\[
|\{x\in A:\dist(x,y)=1\}|\le 2(q-1).
\]

Double counting pairs $(x,y)$ with $x\in A$ and $\dist(x,y)=1$ gives
\[
|A|n(q-1)\le q^n\cdot 2(q-1),
\]
and therefore
\[
|A|\le \frac{2q^n}{n}.
\]

For the lower bound, choose $m$ smallest such that
\[
n\le N:=\frac{q^m-1}{q-1}.
\]
There is a $q$-ary Hamming code $C\subseteq\F_q^N$ with
\[
\dim C=N-m
\qquad\text{and}\qquad
d_{\min}(C)=3;
\]
see \cite[Ch.~7, \S3, Problem~8(a), pp.~193--194]{MS}.

Using the same shortening argument as in the proof of \Cref{lem:bch}, we
obtain a code $C'\subseteq\F_q^n$ with minimum distance at least $3$ and
\[
|C'|\ge q^{n-m}
\]
By the minimality of $m$,
\[
\frac{q^{m-1}-1}{q-1}<n.
\]
It follows that
\[
q^{m-1}<(q-1)n+1,
\]
and hence
\[
q^m=O_q(n).
\]
Therefore,
\[
q^{n-m}
=
\Omega_q\!\left(\frac{q^n}{n}\right).
\]
Since $C'$ contains no pair of distinct vertices at distance $2$, it is an
independent set in $H_q(n,2)$, and therefore triangle-free.  Thus
\[
T_q(n,2)\ge |C'|
=\Omega_q\!\left(\frac{q^n}{n}\right).
\]
\end{proof}

\bigskip

We therefore pose the following problem.
\begin{problem}
    For fixed prime power $q\ge3$ and fixed $r\ge3$, determine $T_q(n,r)$.
\end{problem}

We can also ask the same question for larger cliques.
For
integers $s\ge 3$, define
\[
T_{q,s}(n,r)
:=
\max\Bigl\{
|A|:\ A\subseteq \F_q^n
\text{ and } A \text{ spans no copy of } K_s \text{ in } H_q(n,r)
\Bigr\}.
\]
With this notation, $T_{q,3}(n,r)=T_q(n,r)$.
In the binary setting, the fixed even distance case is also settled for
larger cliques \cite{balogh}.
We close with the following two problems.

\begin{problem}
Fix $s\ge3$ and let $r$ be even. 
Determine $T_{2,s}(n,r)$ when $r=o(n)$ and $r\to\infty$,
and when $r=\alpha n$ with fixed $0<\alpha<2/3$.
\end{problem}

\begin{problem}
Fix a prime power $q\ge3$, an integer $s\ge3$, and a fixed distance
$r\ge1$.  Determine $T_{q,s}(n,r)$.
\end{problem}

\textbf{Acknowledgments.}
The first author is grateful to D\'aniel Varga for introducing her to the problem and for his many helpful suggestions.
The authors also thank Anurag Bishnoi for his helpful comments on the exposition and the proofs.

\end{document}